\documentclass[12pt,a4paper,twoside]{article}

\newcommand{\pageformat}[6]{\setlength{\hoffset}{-1in}
                  \setlength{\voffset}{-1in}
                  \addtolength{\hoffset}{#5}
                            \addtolength{\voffset}{#6}
                            \setlength{\oddsidemargin}{#1}
                            \setlength{\evensidemargin}{#2}
                            \setlength{\textwidth}{\paperwidth}
                  \addtolength{\textwidth}{-\oddsidemargin}
                  \addtolength{\textwidth}{-\evensidemargin}
                  \addtolength{\textwidth}{-\marginparsep}
                  \addtolength{\textwidth}{-\marginparwidth}
                            \setlength{\topmargin}{#3}
                            \setlength{\textheight}{\paperheight}
                  \addtolength{\textheight}{-\topmargin}
                  \addtolength{\textheight}{-\headheight}
                  \addtolength{\textheight}{-\headsep}
                  \addtolength{\textheight}{-\footskip}
                  \addtolength{\textheight}{-#4}}
\pageformat{2cm}{3cm}{24mm}{24mm}{1pt}{0pt}

\usepackage{ifthen}
\newboolean{article}
    \setboolean{article}{true}
\newboolean{report}
\newboolean{book}
\newboolean{letter}
\newboolean{german}
\newboolean{italian}
\newboolean{nobaselinestretch}
\newboolean{nosectionappendix}
\newboolean{oldtoc}
\newboolean{nosectionequn}
\newboolean{notheorem}

\ifthenelse{\boolean{german}}{
    \usepackage{german}}{}

\usepackage[latin1]{inputenc}
\usepackage{hyperref}

\usepackage{amsmath}
\usepackage{amssymb}
\usepackage[mathscr]{eucal}

\ifthenelse{\boolean{notheorem}}{}{
    \usepackage{theorem}}

\ifthenelse{\boolean{nobaselinestretch}}{}{
    \renewcommand{\baselinestretch}{1.25}}

\newenvironment{env}[2]{\begin{#1}#2\end{#1}}{}
    \newcommand{\beq}[1]{\begin{env}{equation}{#1}}
    \newcommand{\beqn}[1]{\begin{env}{equation*}{#1}}
    \newcommand{\bal}[1]{\begin{env}{align}{#1}}
    \newcommand{\baln}[1]{\begin{env}{align*}{#1}}
    \newcommand{\bga}[1]{\begin{env}{gather}{#1}}
    \newcommand{\bgan}[1]{\begin{env}{gather*}{#1}}
    \newcommand{\bflal}[1]{\begin{env}{flalign}{#1}}
    \newcommand{\bflaln}[1]{\begin{env}{flalign*}{#1}}
    \newcommand{\bmu}[1]{\begin{env}{multline}{#1}}
    \newcommand{\bmun}[1]{\begin{env}{multline*}{#1}}
    \newcommand{\bsp}[1]{\begin{env}{split}{#1}}

    \newcommand{\eeq}{\end{env}}
    \newcommand{\eeqn}{\end{env}}
    \newcommand{\eal}{\end{env}}
    \newcommand{\ealn}{\end{env}}
    \newcommand{\ega}{\end{env}}
    \newcommand{\egan}{\end{env}}
    \newcommand{\eflal}{\end{env}}
    \newcommand{\eflaln}{\end{env}}
    \newcommand{\emu}{\end{env}}
    \newcommand{\emun}{\end{env}}
    \newcommand{\esp}{\end{env}}

\newcommand{\lf}{\vspace{2ex}}

\renewcommand{\bf}[1]{\textbf{#1}}
\renewcommand{\it}[1]{\textit{#1}}

\renewcommand{\sf}[1]{\textsf{#1}}

\renewcommand{\tt}[1]{\texttt{#1}}
\newcommand{\hl}[1]{\bf{\it{#1}}}
\newcommand{\mrm}[1]{\mathrm{#1}}

\newcommand{\msf}[1]{\text{\small$\sf{#1}$}}

\newcommand{\cmc}[1]{\mathcal{#1}}
\newcommand{\eus}[1]{\mathscr{#1}}

\newcommand{\bb}[1]{\mathbb{#1}}

\newcommand{\mfootnotesize}[1]{{\setlength{\arraycolsep}{.5ex}\text{\footnotesize$#1$}}}

\newcommand{\mtiny}[1]{{\setlength{\arraycolsep}{.3ex}\text{\tiny$#1$}}}
\newcommand{\nbd}[1]{$#1$\nobreakdash--}
\newcommand{\ol}[1]{\overline{#1}}

\newcommand{\vt}{\vartheta}

\newcommand{\vp}{\varphi}

\newcommand{\norm}[1]{\left\lVert#1\right\rVert}

\newcommand{\bfam}[1]{\bigl(#1\bigr)}
\newcommand{\Bfam}[1]{\Bigl(#1\Bigr)}
\newcommand{\AB}[1]{\langle#1\rangle}
\newcommand{\bAB}[1]{\bigl\langle#1\bigr\rangle}
\newcommand{\BAB}[1]{\Bigl\langle#1\Bigr\rangle}
\newcommand{\CB}[1]{\{#1\}}
\newcommand{\bCB}[1]{\bigl\{#1\bigr\}}
\newcommand{\BCB}[1]{\Bigl\{#1\Bigr\}}

\newcommand{\Matrix}[1]{\begin{pmatrix}#1\end{pmatrix}}

\newcommand{\fMatrix}[1]{\mfootnotesize{\Matrix{#1}}}

\newcommand{\tMatrix}[1]{\mtiny{\Matrix{#1}}}

\newcommand{\set}[2][]{
    \ifthenelse{\equal{#1}{}}{
        \CB{#2}}{
        \CB{#1~|~#2}}}
\newcommand{\bset}[2][]{
    \ifthenelse{\equal{#1}{}}{
        \bCB{#2}}{
        \bCB{#1~|~#2}}}
\newcommand{\Bset}[2][]{
    \ifthenelse{\equal{#1}{}}{
        \BCB{#2}}{
        \BCB{#1~\big|~#2}}}

\DeclareMathOperator{\ls}{\normalfont\msf{span}}

\DeclareMathOperator{\id}{\normalfont\msf{id}}

\DeclareMathOperator{\alg}{\normalfont\msf{alg}}

\newcommand{\dom}{\operatorname{\msf{dom}}}

\newcommand{\C}{\bb{C}}

\newcommand{\E}{\bb{E}}

\newcommand{\N}{\bb{N}}

\newcommand{\R}{\bb{R}}

\newcommand{\cB}{\cmc{B}}
\newcommand{\cC}{\cmc{C}}

\newcommand{\sA}{\eus{A}}
\newcommand{\sB}{\eus{B}}

\newcommand{\sF}{\eus{F}}

\newcommand{\sK}{\eus{K}}

\newcommand{\F}{{\mrm{F}}}

\ifthenelse{\boolean{nosectionequn}}{}{
    \numberwithin{equation}{section}
    }

\ifthenelse{\boolean{article}\or\boolean{letter}\or\boolean{nosectionequn}}{
    \setboolean{nosectionappendix}{true}}{}
\ifthenelse{\boolean{nosectionappendix}}{}{
    \renewcommand{\appendix}{
        \chapter*{\appendixname}
        \addcontentsline{toc}{chapter}{\appendixname}
        \renewcommand{\thesection}{\Alph{section}}
        \setcounter{section}{0}}}
   
\ifthenelse{\boolean{report}\or\boolean{book}}{
    }{}

\ifthenelse{\boolean{notheorem}}{}{
        \newcommand{\mnname}{Mathematical note.}
        \newcommand{\enname}{End of the note.}
        \newcommand{\definame}{Definition.}
        \newcommand{\propname}{Proposition.}
        \newcommand{\lemname}{Lemma.}
        \newcommand{\exname}{Example.}
        \newcommand{\exername}{Exercise.}
        \newcommand{\remname}{Remark.}
        \newcommand{\obname}{Observation.}
        \newcommand{\thmname}{Theorem.}
        \newcommand{\corname}{Corollary.}
        \newcommand{\proofname}{Proof.}
    \ifthenelse{\boolean{german}}{
        \renewcommand{\mnname}{Mathematische Notiz.}
        \renewcommand{\enname}{Ende der Notiz.}
        \renewcommand{\exname}{Beispiel.}
        \renewcommand{\exername}{Übung.}
        \renewcommand{\remname}{Bemerkung.}
        \renewcommand{\obname}{Beobachtung.}
        \renewcommand{\thmname}{Satz.}
        \renewcommand{\corname}{Korollar.}
        \renewcommand{\proofname}{Beweis.}}{}
    \ifthenelse{\boolean{italian}}{
        \renewcommand{\mnname}{Nota matematica.}
        \renewcommand{\enname}{Fina della nota.}
        \renewcommand{\definame}{Definizione.}
        \renewcommand{\propname}{Proposizione.}
        \renewcommand{\exname}{Esempio.}
        \renewcommand{\exername}{Esercizio.}
        \renewcommand{\remname}{Nota.}
        \renewcommand{\obname}{Osservazione.}
        \renewcommand{\thmname}{Teorema.}
        \renewcommand{\corname}{Corollario.}
        \renewcommand{\proofname}{Dimostrazione.}

       \renewcommand{\appendixname}{Appendice}

       }{}
    \theoremheaderfont{\normalfont\bfseries}
    \theoremstyle{change}
        \theorembodyfont{\rmfamily}
            \newtheorem{emp}{}[section]
                \newcommand{\bemp}[1][]{
                    \begin{emp}\hskip-\labelsep\bf{#1}\hskip\labelsep}
                \newcommand{\eemp}{\end{emp}}
\newtheorem{itemp}[emp]{}
                \newcommand{\bitemp}[1][]{
                    \begin{itemp}\hskip-\labelsep\bf{#1}\hskip\labelsep\normalfont\itshape}
                \newcommand{\eitemp}{\end{itemp}}
            \newtheorem{mn}[emp]{\mnname}
                \newcommand{\bnm}{\begin{mn}~\begin{quotation}\renewcommand{\baselinestretch}{1}\small\noindent\ignorespaces}
                \newcommand{\enm}{\end{quotation}\hfill\bf{\enname}\end{mn}}
            \newtheorem{ex}[emp]{\exname}
                \newcommand{\bex}{\begin{ex}}
                \newcommand{\eex}{\end{ex}}
            \newtheorem{exer}[emp]{\exername}
                \newcommand{\bexer}{\begin{exer}}
                \newcommand{\eexer}{\end{exer}}
            \newtheorem{defi}[emp]{\definame}
                \newcommand{\bdefi}{\begin{defi}}
                \newcommand{\edefi}{\end{defi}}
            \newtheorem{rem}[emp]{\remname}
                \newcommand{\brem}{\begin{rem}}
                \newcommand{\erem}{\end{rem}}
            \newtheorem{ob}[emp]{\obname}
                \newcommand{\bob}{\begin{ob}}
                \newcommand{\eob}{\end{ob}}
        \theorembodyfont{\normalfont\itshape}
            \newtheorem{thm}[emp]{\thmname}
                \newcommand{\bthm}{\begin{thm}}
                \newcommand{\ethm}{\end{thm}}
            \newtheorem{prop}[emp]{\propname}
                \newcommand{\bprop}{\begin{prop}}
                \newcommand{\eprop}{\end{prop}}
            \newtheorem{cor}[emp]{\corname}
                \newcommand{\bcor}{\begin{cor}}
                \newcommand{\ecor}{\end{cor}}
            \newtheorem{lem}[emp]{\lemname}
                \newcommand{\blem}{\begin{lem}}
                \newcommand{\elem}{\end{lem}}
\newenvironment{empn}[1]{\lf\noindent\bf{#1}\ignorespaces\hskip\labelsep}{\lf}
		\newcommand{\bempn}[1]{\begin{empn}{#1}}
		\newcommand{\eempn}{\end{empn}}
		\newcommand{\bitempn}[1]{\begin{empn}{#1}\normalfont\itshape}
		\newcommand{\eitempn}{\end{empn}}
                \newcommand{\bnmn}{\begin{empn}{\mnname}~\begin{quotation}\renewcommand{\baselinestretch}{1}\small\noindent\ignorespaces}
                \newcommand{\enmn}{\end{quotation}\hfill\bf{\enname}\end{empn}}
		\newcommand{\bexn}{\begin{empn}{\exname}}
		\newcommand{\eexn}{\end{empn}}
		\newcommand{\bexern}{\begin{empn}{\exername}}
		\newcommand{\eexern}{\end{empn}}
		\newcommand{\bdefin}{\begin{empn}{\definame}}
		\newcommand{\edefin}{\end{empn}}
		\newcommand{\bremn}{\begin{empn}{\remname}}
		\newcommand{\eremn}{\end{empn}}
		\newcommand{\bobn}{\begin{empn}{\obname}}
		\newcommand{\eobn}{\end{empn}}

\newcommand{\qedsymbol}{~\rule[-0.35mm]{2mm}{2mm}}
    \newcounter{proof}[emp]
    \newenvironment{Proof}[1]{
        \vspace{1ex}
        \renewcommand{\item}[1][\stepcounter{proof}(\roman{proof})]%
            {##1\hskip\labelsep}
        \noindent\textsc{#1\hskip\labelsep}}{
        \nolinebreak\qedsymbol}
    \newcommand{\proof}[1][\proofname]{
        \begin{Proof}{#1}\ignorespaces}
    \newcommand{\qed}{\end{Proof}}
    \newcommand{\noqed}{
        \renewcommand{\qedsymbol}{}
        \end{Proof}}}
    \ifthenelse{\boolean{italian}}{
        \renewcommand{\proofname}{Dimostrazione.}}{}

\usepackage[varg]{txfonts}

\begin{document}




\title{Generators of Dynamical Systems on Hilbert Modules}
\author{Gholamreza Abbaspour Tabadkan\thanks{GAT is supported by a grant of Ferdowsi University, Mashhad.}~~ and~ Michael Skeide\thanks{MS is supported by research funds of the University of Molise and the Italian MIUR (PRIN 2005).}}

\date{October 2006}

{
\renewcommand{\baselinestretch}{1}
\maketitle



\begin{abstract}\noindent
We characterize the generators of dynamical systems on Hilbert modules as those generators of one-parameter groups of Banach space isometries which are ternary derivations. We investigate in how far a similar condition can be expressed in terms of generalized derivations.
\end{abstract}


}


%

\section{Introduction} \label{intro}

Let $E$ be a Hilbert module over a \nbd{C^*}algebra $\cB$. A \hl{generalized unitary} on $E$ is a surjection $u$ on $E$ that fulfills
\beqn{\tag{GU}\label{genu}
\AB{ux,uy}
~=~
\vp(\AB{x,y}),
~~~~~~
x,y\in E
}\eeq
for some automorphism $\vp$ of $\cB$. We will also say $u$ is a \hl{\nbd{\vp}unitary}. A \hl{generalized derivation} of $E$ is a densely defined linear map $\delta\colon E\supset\dom(\delta)\rightarrow E$ that fulfills
\beqn{\tag{GD}\label{gend}
\delta(xb)
~=~
\delta(x)b+xd(b),
~~~~~~
x\in\dom(\delta),b\in\dom(d)
}\eeqn
for some derivation $d\colon\cB\supset\dom(d)\rightarrow\cB$ of $\cB$, in such a way that $\dom(\delta)$ is a right $\dom(d)$--mod\-ule. We will also say $\delta$ is a \hl{\nbd{d}derivation}. A \hl{dynamical system} on a Hilbert \nbd{\cB}module $E$ is a strongly continuous one-parameter group $u=\bfam{u_t}_{t\in\R}$ of generalized unitaries. Abbaspour, Moslehian and Niknam \cite{AMN05} showed that the generator of a dynamical system is a generalized derivation. However, even if a closed and densely defined map $\delta$ on $E$ is the generator of a group of Banach space isometries, then for that this group forms a dynamical system, it is not sufficient that $\delta$ be a generalized derivation.

It is the scope of these notes to find a better algebraic condition. This condition will be in terms of \it{ternary} maps.

A \hl{ternary automorphism} of $E$ is a bijection $u$ on $E$ that fulfills
\beqn{\tag{TU}\label{teru}
u(x\AB{y,z})
~=~
(ux)\AB{uy,uz},
~~~~~~
x,y,z\in E.
}\eeqn
In Section \ref{isoSEC} we show that the generalized isometries from a \it{full} Hilbert \nbd{\cB}module $E$ to a Hilbert \nbd{\cC}module $F$ are exactly the ternary homomorphims. As a special case this includes the statement that the generalized unitaries on a full Hilbert \nbd{\cB}module are exactly its ternary automorphisms. This frees the discussion from worrying about existence of an automorphism $\vp$ of $\cB$. In fact, the main problem in the proof is to show that existence of such an automorphism is automatic. Consequently, the dynamical systems on a full Hilbert module $E$ are exactly the strongly continuous one-parameter groups of ternary automorphisms.

A \hl{ternary derivation} of $E$ is a densely defined linear map 
$\delta\colon E\supset\dom(\delta)\rightarrow E$ that fulfills
\beqn{\tag{TD}\label{terd}
\delta(x\AB{y,z})
~=~
\delta(x)\AB{y,z}+x\AB{\delta(y),z}+x\AB{y,\delta(z)}
~~~~~~
x,y,z\in E
}\eeqn
where $\dom(\delta)~\bAB{\,\dom(\delta)\,,\,\dom(\delta)\,}\subset\dom(\delta)$, that is, $\dom(\delta)$ is invariant under the \hl{ternary product} $(x,y,z)\mapsto x\AB{y,z}$. In Section \ref{derSEC} we show that every ternary derivation of a full Hilbert module is a generalized derivation, while the converse fails. Generators of dynamical systems are always ternary derivations. We show also a sort of converse: If a linear densely defined map on $E$ is the generator of a \hl{\nbd{C_0}group} (that is, a strongly continuous one-parameter group of Banach space isometries) on $E$, then this group is a dynamical system, if and only if $\delta$ is a ternary derivation. This reduces the problem of characterizing the generators to the well-known general analytic criteria based on Hille-Yosida theory that state when $\delta$ is the generator of a \nbd{C_0}group, and the purely algebraic question whether $\delta$ is a ternary derivation. We see that we have a satisfactory description of generators of dynamical systems on Hilbert modules in terms of ternary derivations, while the larger part of Section \ref{derSEC} illustrates that similar statements in terms of generalized derivations are possible only under rather hard analytical hypothesis.

We note, too, that the condition that the Hilbert \nbd{\cB}module $E$ be full is not critical as long as we speak about ternary maps. Restrictions that arise in the case of generalized unitaries \it{on} a Hilbert module $E$ when $\cB$ is not chosen minimal have been analyzed in Skeide \cite{Ske05p2}.

\lf\noindent
\bf{Acknowledgements.~}
The results in Section \ref{isoSEC} were included as a part of the first authors PhD-thesis \cite{Abb06}. Most results have been obtained during a six months visit of the first author at the Dipartimento S.E.G.e S.\, financed by a grant from Ferdowsi University, Mashhad. Both authors acknowledge the support by research funds of the Italian MIUR and University of Molise.

\lf\noindent
\bf{Conventions and notations.~}
A \hl{pre-Hilbert \nbd{\cB}module} is a right module $E$ over a (pre-)\nbd{C^*}al\-ge\-bra, with a sesquilinear inner product $\AB{\bullet,\bullet}\colon E\times E\rightarrow\cB$ that satisfies $\AB{x,yb}=\AB{x,y}b$ $(x,y\in E;b\in\cB)$, $\AB{x,x}\ge0$ and $\AB{x,x}=0$ $\Rightarrow$ $x=0$ $(x\in E)$. A \hl{Hilbert \nbd{\cB}module} is a pre-Hilbert \nbd{\cB}module that is complete in the \hl{norm} $\norm{x}:=\sqrt{\AB{x,x}}$. A pre-Hilbert \nbd{\cB}module $E$ is \hl{full}, if the the \hl{range ideal} $\cB_E:=\ls\AB{E,E}$ is dense in $\cB$.

By $\sB^a(E)$ and $\sK(E)$ we denote the \nbd{C^*}algebras of all \hl{adjointable} operators and of all \hl{compact} operators, respectively, on $E$, where $\sK(E)$ is the completion of the pre-\nbd{C^*}algebra $\sF(E)$ of \hl{finite-rank} operators which is spanned by the \hl{rank-one} operators $xy^*\colon z\mapsto x\AB{y,z}$.

\section{Generalized isometries \it{versus} ternary homomorphisms}\label{isoSEC}

Unitaries on or between Hilbert modules are inner product preserving surjections. For isometries, surjectivity is missing. For generalized unitaries on a Hilbert module in \eqref{genu} the condition that the surjection preserves inner products is modified to that it preserves inner products up to a fixed automorphism of the algebra. When the unitary is between different Hilbert modules, it is not even necessary that these are modules over the same algebra. In this section we investigate generalized isometries between Hilbert modules.

Let $E$ be a Hilbert \nbd{\cB}module and let $F$ be a Hilbert \nbd{\cC}module. A \hl{generalized isometry} from $E$ to $F$ is a map $u\colon E\rightarrow F$ that fulfills \eqref{genu} for some homomorphism $\vp\colon\cB\rightarrow\cC$. We will also say $u$ is a \hl{\nbd{\vp}isometry}. Calculating the norm of $ux+(uy)\vp(b)-u(x+yb)$, we find that every \nbd{\vp}isometry $u$ is \hl{\nbd{\vp}linear}, that is, $ux+(uy)\vp(b)=u(x+yb)$ $(x,y\in E;b\in\cB)$. Inserting scalar multiples of an approximate unit, we see that \nbd{\vp}linearity implies \nbd{\C}linearity. Obviously, a \nbd{\vp}isometry has norm $1$, unless $\vp\upharpoonright\cB_E$ is $0$.

A \hl{ternary homomorphism} from $E$ to $F$ is a map $u\colon E\rightarrow F$ that fulfills \eqref{teru}. Obviously, every \nbd{\vp}isometry is a ternary homomorphism. It is our scope to show that, at least if $E$ is full, then every ternary homomorphism is also a \nbd{\vp}isometry.

\bthm\label{gen=ter}
For a map $u$ from a full Hilbert \nbd{\cB}module $E$ to a Hilbert \nbd{\cC}module $F$ the following statements are equivalent:
\begin{enumerate}
\item
$u$ is a generalized isometry.

\item
$u$ is a ternary homomorphism.
\end{enumerate}
\ethm

\proof
Given a ternary homomorphism $u$ from a full Hilbert \nbd{\cB}module $E$ to a Hilbert \nbd{\cC}mod\-ule $F$, it is our job to find a homomorphism $\vp\colon\cB\rightarrow\cC$ such that $u$ fulfills \eqref{genu}. First, we observe that for full $E$ such a homomorphism is determined uniquely by \eqref{genu}. The attempt to define the homomorphism $\vp$ first on the pre-\nbd{C^*}algebra $\cB_E$ by $\AB{x,y}\mapsto\AB{ux,uy}$ and then to show that it is bounded by appealing to Muhly, Skeide and Solel \cite[Corollary 1.20]{MSS06} has been put into practise in \cite{Abb06} under the assumption that $u$ is linear. Here we follow a different road.

Let us observe that if a suitable $\vp$ exists, then $u$ must be \nbd{\vp}linear. So, necessarily we must have $(ux)\vp(b)=u(xb)$ for all $x\in E$. We use this property to define a left action $\vp(b)$ of $b\in\cB$ on the pre-\nbd{C^*}algebra $\cC_{uE}:=\ls\AB{uE,uE}$ considered pre-Hilbert module over itself in the usual way, that is, with inner product $\AB{c,c'}=c^*c'$ and right action simply by multiplication. We put
\beqn{
\vp(b)\AB{ux,uy}
~:=~
\AB{u(xb^*),uy}
}\eeqn
and we must show, in a first step, that this well-defines a homomorphism into $\sB^a(\ol{\cC_{uE}})$. As, clearly, $\vp(b)\vp(b')\AB{ux,uy}=\vp(bb')\AB{ux,uy}$ (so that, once well-defined, $\vp$ is multiplicative), it is sufficient to show that $\vp(b^*)$ is a formal adjoint of $\vp(b)$ on the spanning subset of elements of the form $\AB{ux,uy}$. From this follow both that $\vp(b)$ is well-defined and that $\vp(b^*)=\vp(b)^*$. We start by observing that
\beqn{
\bAB{c,\AB{ux,uy}}
~=~
c^*\AB{ux,uy}
~=~
\AB{(ux)c,uy}
}\eeqn
for all $c\in\cC_{uE}$. Using this two times, we find
\bmun{
\BAB{\AB{ux,uy},\vp(b)\AB{ux',uy'}}
~=~
\BAB{\AB{ux,uy},\AB{u(x'b^*),uy'}}
~=~
\bAB{u(x'b^*)\AB{ux,uy},uy'}
\\
~=~
\bAB{u(x'b^*\AB{x,y}),uy'}
~=~
\bAB{u(x'\AB{xb,y}),uy'}
\\
~=~
\bAB{(ux')\AB{u(xb),uy}),uy'}
~=~
\BAB{\AB{u(xb),uy},\AB{ux',uy'}}
~=~
\BAB{\vp(b^*)\AB{ux,uy},\AB{ux',uy'}}.
}\emun
Like every homomorphism from a \nbd{C^*}algebra into the adjointable operators on a pre-Hilbert module, $\vp$ maps into the bounded operators, and like every homomorphism from a \nbd{C^*}algebra into a pre-\nbd{C^*}algebra, $\vp$ is a contraction.

Next we observe that $\vp(\AB{x,y})$ acts on the element $\AB{ux',uy'}$ of $\cC_{uE}$ simply by multiplication from the left with the element $\AB{ux,uy}$. The subalgebra $\vp(\cB_E)$ of $\sB^a(\ol{\cC_{uE}})$ is nothing but $\cC_{uE}$, which, of course, is faithfully contained in $\sB^a(\ol{\cC_{uE}})$. In other words, $\vp$ is the unique continuous extension from $\cB_E$ to $\cB=\ol{\cB_E}$ of $\vp\upharpoonright\cB_E$ and, therefore, maps into $\ol{\cC_{uE}}\subset\cC$. Clearly, $\vp$ turns $u$ into a \nbd{\vp}isometry.\qed

\bcor
Every ternary homomorphism is linear and contractive.
\ecor

\proof
The only thing that remains is to remark that if $E$ is not full, then we may simply turn $E$ into a full Hilbert module by restricting to $\ol{\cB_E}$.\qed

\bob
Note that a ternary homomorphism is injective, if and only if the homomorphism $\vp\colon\ol{\cB_E}\rightarrow\cC$ that turns it into a \nbd{\vp}isometry is injective. (Every surjective but noninjective endomorphism of $\cB$ is an example for a surjective ternary homomorphism that is not injective.) This shows, in particular, that the $\vp$ induced by a ternary automorphism on a full Hilbert module is itself an automorphism.
\eob

\bcor\label{dscor}
The group of generalized unitaries on a full Hilbert module $E$ coincides with the group of ternary automorphisms of $E$. Therefore, the dynamical systems on a Hilbert module $E$ are exactly the \nbd{C_0}groups of ternary automorphisms.
\ecor

\brem\label{terrem}
By the construction in the proof of Theorem \ref{gen=ter} every \nbd{C_0}group $u=\bfam{u_t}_{t\in\R}$ of ternary automorphisms of a Hilbert \nbd{\cB}module comes along with a (unique) family of automorphisms $\vp_t$ of $\ol{\cB_E}$ and, obviously, the $\vp_t$ form a \nbd{C^*}dynamical system. These automorphisms $\vp_t$ do, in general, not necessarily extend to automorphisms of $\cB$; see \cite{Ske05p2}. Therefore, for not necessarily full $E$ there are, in general, more groups of ternary automorphisms than groups of generalized unitaries. In the general case, it seems, therefore, advisable to define a dynamical system on a Hilbert module as a \nbd{C_0}group of ternary automorphisms.
\erem

\brem\label{linkrem}
By \cite[Observation 1.4]{Ske05p2} (for instance) we know that every surjective \nbd{\vp}iso\-me\-try from a Hilbert \nbd{\cB}module $E$ to a Hilbert \nbd{\cC}module $F$ extends to a homomorphism between the \hl{extended linking algebras}
\beqn{
\Phi
~=~
\fMatrix{\vp&u^*\\u&\vt}
~\colon~
\fMatrix{\cB&E^*\\E&\sB^a(E)}
~\longrightarrow~
\fMatrix{\cC&F^*\\F&\sB^a(F)}
}\eeqn
that restricts to a homomorphism between the usual \hl{linking algebras} $\tMatrix{\cB&E^*\\E&\sK(E)}\rightarrow\tMatrix{\cC&F^*\\F&\sK(F)}$. (Here $u^*(x^*):=(ux)^*$, while $\vt(a)$ acts on $y=ux$ in the only possible way, namely, $\vt(a)(ux)=u(ax)$. Well-definedness of $\vt(a)$ follows in a way paralleling the proof of well-definedness of $\vp$ in the proof of Theorem \ref{gen=ter}.) Therefore, generalized isometries and, consequently, also ternary homomorphisms are even completely contractive. (One may obtain this result also as in \cite{Abb06}, by showing that every inflation $u^n$ of $u$ is a \nbd{\vp^n}isometry from $M_n(E)$ to $M_n(F)$ and, therefore, a contraction.) This improves on a result on ternary homomorphisms of \nbd{C^*}algebras by Bracic and Moslehian \cite{BrMo06p}.

A ternary homomorphism $\eta$  from $E$ into the Hilbert \nbd{\sB(G)}module $\sB(G,H)$ for two Hilbert spaces $G$ and $H$ is what has been called a \hl{representation} of $E$ from $G$ to $H$ in Skeide \cite{Ske00b}. The preceding discussion improves also on \cite[Theorem A.4]{Ske00b} where the extendibility of $\eta$ to a representation of the linking algebra has been shown under the explicit hypothesis that $\eta$ be completely bounded. Now we see that this hypothesis is fulfilled automatically.
\erem

\section{Generalized derivations \it{versus} ternary derivations}\label{derSEC}

It is easy to see that the generator $\delta$ of a dynamical system $u=\bfam{u_t}_{t\in\R}$ on a full Hilbert \nbd{\cB}module $E$ is a generalized derivation; see \cite{AMN05} and cf.\ also Corollary \ref{AMNcor}. A possible choice for the derivation $d$ in \eqref{gend} is simply the generator of the \nbd{C^*}dynamical system $\vp=\bfam{\vp_t}_{t\in\R}$ \hl{associated} with $u$, that is determined uniquely by the requirement that every $u_t$ is a \nbd{\vp_t}unitary; see Remark \ref{terrem}. But, even if we know that the generator $\delta$ of \nbd{C_0}group on $E$ is a \nbd{d}derivation, then it is not possible to conclude that $\delta$ generates a dynamical system without making further analytical assumptions about $d$ and algebraic assumptions about the domain of $d$, see Theorem \ref{dclosethm}. These algebraic conditions are \it{relative} to $\delta$, that is, they cannot be formulated intrinsically in terms of the derivation $d$ of $\cB$ alone, but depend explicitly on $\delta$. On the other hand,  it is easy to formulate these conditions intrinsically in terms of $\delta$ alone: $\delta$ must be a ternary derivation.

We study, first, the intrinsic description of the generators of dynamical systems on Hilbert modules as ternary derivations (Theorem \ref{genterthm}). Then, we explain the relationship between ternary derivations and generalized derivations. We will see that there is a particular derivation $d_\delta$ (Theorem \ref{tergenthm}) that allows to formulate Theorem \ref{genterthm} in terms of generalized derivations (Theorem \ref{gengenthm}). In Theorem \ref{sumthm} we summarize all criteria and add one more in terms of the linking algebra.

\bthm\label{genterthm}
Let $u=\bfam{u_t}_{t\in\R}$ be \nbd{C_0}group on a Hilbert \nbd{\cB}module $E$. Then $u$ is a dynamical system if and only if its generator $\delta$ is a ternary derivation.
\ethm

\proof
Recall that the generator of a \nbd{C_0}group $u$ is defined as
\beqn{
\delta(x)
~:=~
\lim_{t\to0}\frac{u_tx-x}{t}
}\eeqn
for all $x$ for which the limit exists. Further, recall that this domain $\dom(\delta)$ contains a  dense \hl{core} of \hl{entire analytic vectors}. That means, the subspace $\sA(\delta)\subset\bigcap_{n\in\N}\dom(\delta^n)$ that consists of all vectors $x$ for which for all $t\in\R$ the series
\beqn{
\sum_{n=0}^\infty\frac{t^n\delta^n}{n!}x
}\eeqn
converges absolutely to the limit $u_tx$ is dense in $E$ and $\delta$ is the closure of $\delta\upharpoonright\sA(\delta)$. 

Suppose $\delta$ is the generator of a dynamical system $u$. Let $x,y,z\in\dom(\delta)$. By Corollary \ref{dscor} all $u_t$ are ternary automorphisms, so that
\bmun{
\frac{u_t(x\AB{y,z})-x\AB{y,z}}{t}
~=~
\frac{(u_tx)\AB{u_ty,u_tz}-x\AB{y,z}}{t}
\\
~=~
\frac{u_tx-x}{t}\AB{u_ty,u_tz}+x\BAB{\frac{u_ty-y}{t},u_tz}+x\BAB{y,\frac{u_tz-z}{t}}.
}\emun
As all $u_t$ are contractions, the families $u_tx$ and $u_ty$ are bounded uniformly. So the limit of the preceding expression exists and is equal to $\delta(x)\AB{y,z}+x\AB{\delta(y),z}+x\AB{y,\delta(z)}$. This shows both that $x\AB{y,z}\in\dom(\delta)$ and that $\delta$ is a ternary derivation.

Conversely, suppose that $\delta$ is a ternary derivation, and choose entire analytic elements $x,y,z\in\sA(\delta)$. By a routine induction we show the ternary generalized \it{Leibniz rule}
\beqn{
\delta^n(x\AB{y,z})
~=~
\sum_{\substack{k,\ell,m\in\N_0\\k+\ell++m=n}}\frac{n!}{k!\ell!m!}\delta^k(x)\AB{\delta^\ell(y),\delta^m(z)}.
}\eeqn
From this, one easily concludes that $x\AB{y,z}$ is also in $\sA(\delta)$ and that $u_t$ fulfills \eqref{teru} on the dense subset $\sA(\delta)$. By contractivity of $u_t$, this extends to all of $E$ so that $u_t$ is a ternary automorphism.\qed

\brem
For general results about \nbd{C_0}groups we refer to Bratteli and Robinson \cite{BrRo87}. In particular, the problem to decide, whether a linear densely defined map is the generator of \nbd{C_0}group, we leave entirely to the comprehensive treatment in \cite{BrRo87}. But, once we have such a generator, we see that the problem whether the generated group is a dynamical system, is equivalent to the question whether the generator is a ternary derivation. Thus, we have a complete separation into the general analytic criteria of the Banach space theory that determine when $\delta$ is a generator (which we do not treat here) and the completely algebraic criterion in Theorem \ref{genterthm}.
\erem

Theorem \ref{genterthm}, in principle, completely settles the problem to characterize the generators of dynamical systems on Hilbert modules. Fullness, is not at all a critical assumption, because if necessary we may always make $\cB$ smaller. In the remainder of this section we deal with the problem to find similarly useful criteria in terms of generalized derivations. We start by establishing a connection between ternary derivations and a special sort of generalized derivations on the algebraic level. However, a full correspondence we will obtain only under the assumption that the derivations in question generate \nbd{C_0}groups. On the level of derivations the assumption of fullness becomes much more vital, as we do not see a possibility to show that the derivation of $\cB$ that turns a map $\delta$ into a \nbd{d}derivation restricts to a derivation of $\ol{\cB_E}$. The following uniqueness result, depending essentially on fullness, is crucial for all other statements that follow.

\bprop\label{duniprop}
Let $\delta\colon E\supset\dom(\delta)\rightarrow E$ be a densely defined linear map on a full Hilbert \nbd{\cB}module $E$. Then for every dense subalgebra $\cB_0$ of $\cB$, there is at most one derivation $d$ of $\cB$ with domain $\dom(d)=\cB_0$ that turns $\delta$ into a \nbd{d}derivation.
\eprop

\bcor
If $d_1,d_2$ are derivations of $\cB$ and if $\delta$ is a \nbd{d_1}derivation \hl{and} a \nbd{d_2}derivation of a full Hilbert \nbd{\cB}module $E$, then
\beqn{
d_1
~\subset~
d_2
~~~~~~\Longleftrightarrow~~~~~~
\dom(d_1)
~\subset~
\dom(d_2).
}\eeqn
\ecor

\proof[Proof of Proposition \ref{duniprop}.~]
If $\delta$ is a \nbd{d}derivation, then we have
\beqn{
xd(b)
~=~
\delta(xb)-\delta(x)b
}\eeqn
for all $x\in\dom(\delta)$ and all $b\in\dom(d)$. Since $E$ is full and $\dom(\delta)$ is dense in $E$, the preceding equation determines $d(b)\in\cB$ uniquely.\qed

\bthm\label{tergenthm}
Every ternary derivation $\delta$ of a full Hilbert \nbd{\cB}module $E$ is also a generalized derivation. More precisely, there is a unique derivation $d_\delta$ of $\cB$ on the dense domain $\dom(d_\delta):=\ls\AB{\dom(\delta),\dom(\delta)}$ that fulfills
\beq{\label{ddef}
d_\delta(\AB{x,y})
~=~
\AB{\delta(x),y}+\AB{x,\delta(y)}.
}\eeq
$d_\delta$ turns $\delta$ into a \nbd{d_\delta}deri\-vation. Moreover, $d_\delta$ is a \nbd{*}derivation.
\ethm

\proof
Suppose we have a derivation $d_\delta$ on the given domain, that turns $\delta$ into a \nbd{d_\delta}derivation. Then (following the proof of Proposition \ref{duniprop}) for the uniquely determined values of $d_\delta(\AB{x,y})$ we find
\beq{\label{dcheck}
x\,d_\delta\bfam{\AB{y,z}}
~=~
\delta(x\AB{y,z})-\delta(x)\AB{y,z}
~=~
x\,\Bfam{\AB{\delta(y),z}+\AB{y,\delta(z)}}
}\eeq
for all $x,y,z\in\dom(\delta)$. We see that if a suitable derivation $d_\delta$ exists, then it must fulfill \eqref{ddef}. In particular, $d_\delta$ is necessarily a \nbd{*}derivation. So the only remaining questions are, firstly, whether \eqref{ddef} always well-defines a linear map
\beqn{
d_\delta
\colon
\ls\AB{\dom(\delta),\dom(\delta)}
~\longrightarrow
~
\cB,
}\eeqn
and, secondly, whether this map is a (necessarily \nbd{*}) derivation. For the first question, suppose $y_i,z_i$ are finitely many elements of $\dom(\delta)$ fulfilling $\sum_i\AB{y_i,z_i}=0$. Then
\beqn{
x\,\Bfam{\sum_i\AB{\delta(y_i),z_i}+\AB{y_i,\delta(z_i)}}
~=~
\delta\Bfam{x\sum_i\AB{y_i,z_i}}-\delta(x)\sum_i\AB{y_i,z_i}
~=~
0
}\eeqn
for all $x\in\dom(\delta)$, so that $d_\delta$ is, indeed, well-defined. For the second question, let us compute the inner product of an element $w\in\dom(\delta)$ with \eqref{dcheck} and the adjoint of the resulting equation. Using this, we find
\bmun{
d_\delta\bfam{\AB{w,x}}\,\AB{y,z}+\AB{w,x}\,d_\delta\bfam{\AB{y,z}}
~=~
\Bfam{\AB{\delta(w),x}+\AB{w,\delta(x)}}\AB{y,z}+\AB{w,x}\,\Bfam{\AB{\delta(y),z}+\AB{y,\delta(z)}}
\\
~=~
\AB{\delta(w),x}\AB{y,z}+\BAB{w\,,\,\delta(x)\AB{y,z}+x\AB{\delta(y),z}+x\AB{y,\delta(z)}}
\\
~=~
\bAB{\delta(w),x\AB{y,z}}+\bAB{w\,,\,\delta(x\AB{y,z})}
~=~
d_\delta\Bfam{\bAB{w,x\AB{y,z}}}
~=~
d_\delta\bfam{\,\AB{w,x}\AB{y,z}\,}.
}\emun
By linearity this extends to $d_\delta(b)b'+bd_\delta(b')=d_\delta(bb')$ for all $b,b'\in\ls\AB{\dom(\delta),\dom(\delta)}$.\qed

\bitemp[Corollary \cite{AMN05}.~]\label{AMNcor}
Every generator of a dynamical system on a full Hilbert module is a generalized derivation.
\eitemp

\bcor\label{linkcor}
Every ternary derivation of a full Hilbert module extends as a \nbd{*}derivation to the linking algebra.
\ecor

\proof
Suppose $\delta$ is ternary derivation of the full Hilbert \nbd{\cB}module. By Theorem \ref{tergenthm} this determines the \nbd{*}derivation $d_\delta$ of $\cB$ which is the candidate for how to extend $\delta$ to the corner $\cB$ of the linking algebra. To find the extension to $\sK(E)$, we observe that $E^*$ is a full Hilbert \nbd{\sK(E)}module, the \hl{dual} module of $E$, with inner product $\AB{x^*,y^*}:=xy^*\in\sK(E)$ and the right action $x^*a:=(a^*x)$ of elements $a\in\sK(E)$ (or even in $\sB^a(E)$). Of course, $\delta^*(x^*):=\delta(x)^*$ defines a ternary derivation of $E^*$ with domain $\dom(\delta^*):=\dom(\delta)^*$, and by Theorem \ref{tergenthm} there is a unique \nbd{*}derivation $d_{\delta^*}$ of $\sK(E)$ defined on the domain $\dom(d_{\delta^*})=\ls\bfam{\dom(\delta)\dom(\delta)^*}$, turning $\delta^*$ into a \nbd{d_{\delta^*}}derivation. It is routine to check that $\tMatrix{b&y^*\\x&a}\mapsto\tMatrix{d_\delta(b)&\delta^*(y^*)\\\delta(x)&d_{\delta^*}(a)}$ defines a \nbd{*}derivation of the linking algebra.\qed

\lf
If $\delta$ is a ternary derivation, the derivation $d_\delta$ plays a distinguished role as it is related more directly to questions of closability than any other derivation $d$ that turns $\delta$ into a \nbd{d}derivation. The following Proposition \ref{terdprop} settles some of these closability questions in the setting of general derivations, while in Theorems \ref{dclosethm} and \ref{gengenthm} the assumption that the maps are generators of \nbd{C_0}groups is crucial. The following task, needed in the proofs of Proposition \ref{terdprop}\eqref{pb} and of Lemma \ref{philinlem}, is so useful that we prefer to formulate it separately.

\blem\label{nlem}
Suppose the elements of the Hilbert \nbd{\cB}module $E$ \hl{separate the points} of $\cB$, that is, $xb=0$ for all $x\in E$ implies $b=0$. (For instance, suppose $E$ is full.) Then
\beqn{
\norm{b}
~=~
\sup_{\norm{x}\le1}\norm{xb}.
}\eeqn
\elem

\proof
By setting $bx^*:=(xb^*)^*$ we define a representation of $\cB$ by adjointable (and, therefore, bounded; see the proof of Theorem \ref{gen=ter}) operators on the dual module $E^*$ (see the proof of Corollary \ref{linkcor}). By hypothesis, this representation is faithful and, therefore, isometric. In other words, the operator norm of the action of $b\in\cB$ as operator on $E^*$ coincides with the norm of $b$ as element of $\cB$. Observing that $\norm{xb}=\norm{(xb)^*}$, this is precisely the statement of the lemma.\qed

\bprop\label{terdprop}
Let $E$ be a full Hilbert \nbd{\cB}module.
\begin{enumerate}
\item\label{p1}
Let $\delta$ be a \nbd{d}derivation of $E$.
\begin{enumerate}
\item\label{pa}
If $\delta$ is closable, then so is $d$.

\item\label{pb}
If $\delta$ is bounded, then so is $d$.
\end{enumerate}

\item\label{p2}
Let $\delta$ be a ternary derivation of $E$. Then $\delta$ is closable, if and only if $d_\delta$ is closable.

\item\label{p3}
Let $\delta$ be a ternary derivation \hl{and} a \nbd{d}derivation of $E$. If $d_\delta$ is closable, then so is $d$.
\end{enumerate}
\eprop

\proof
\item[\eqref{pa}]
Suppose that $\delta$ is a closable \nbd{d}derivation. Let $b_n\to0$ be a sequence in $\dom(d)$ such that $d(b_n)\to b\in\cB$. Then for every $x\in\dom(\delta)$ we find
\beqn{
\delta(xb_n)
~=~
\delta(x)b_n+xd(b_n)
~\longrightarrow~
0+xb.
}\eeqn
As $xb_n\to0$ and $\delta$ is closable, it follows that $\delta(xb_n)\to0$, so that $xb=0$ for all $x\in\dom(\delta)$. As $E$ is full, this implies $b=0$. So, $d$ is closable.

\item[\eqref{pb}]
Suppose that $\delta$ is a bounded \nbd{d}derivation. By Lemma \ref{nlem}, for every $b\in\dom(d)$ we find an $x$ in the unit ball of $E$ such that $\norm{xd(b)}\ge\frac{1}{2}\norm{d(b)}$. So, $\norm{d(b)}\le2\norm{xd(b)}\le2\bfam{\norm{\delta(xb)}+\norm{\delta(x)b}}\le4\norm{\delta}\norm{b}$.

\item[\eqref{p2}]
Suppose now that $\delta$ is a ternary derivation such that $d_\delta$ is closable. Let $x_n\to0$ be a sequence in $\dom(\delta)$ such that $\delta(x_n)\to x\in E$. Then for every $y\in\dom(\delta)$ we find
\beqn{
d_\delta(\AB{y,x_n})
~=~
\AB{y,\delta(x_n)}+\AB{\delta(y),x_n}
~\longrightarrow~
\AB{y,x}+0.
}\eeqn
As $\AB{y,x_n}\to0$ and $d_\delta$ is closeable, it follows that $d_\delta(\AB{y,x_n})\to0$, so that $\AB{y,x}=0$ for all $y\in\dom(\delta)$ and, therefore, $x=0$. So, $\delta$ is closable. If $E$ is full, then, by Part \eqref{pa}, also the converse is true.

\item[\eqref{p3}]
If $d_\delta$ is closable, then, by \eqref{p2}, $\delta$ is closable so that, by \eqref{pa}, $d$ is closable.\qed

\brem
Boundedness of $d$ is not sufficient for boundedness of $\delta$. In fact, every generator of a unitary \nbd{C_0}group on a Hilbert module that is not uniformly continuous is an unbounded ternary derivation and a \nbd{0}derivation for the trivial derivation $0\colon b\mapsto 0$.
\erem

\bob\label{d'ob}
If, in \eqref{p3}, $\dom(d)$ does not contain $\dom(d_\delta)$, then we may easily replace $d$ by the derivation $d'$ defined on $\alg^*\bfam{\dom(d_\delta),\dom(d)}$, the \nbd{*}algebra generated by $\dom(d_\delta)$ and $\dom(d)$, that is determined uniquely (see Proposition \ref{duniprop}!) by the requirement that $\delta$ be a \nbd{d'}derivation. (If such a $d'$ exists, then, again by Proposition \ref{duniprop}, this implies also that $d'$ is the unique extension as a derivation of $d$ and $d_\delta$ to the new domain.) Let us first define $d'$ on the domain $\dom(d)\cup\dom(d_\delta)$ as $d'(b):=d(b)$ for $b\in\dom(d)$ and $d'(b):=d_\delta(b)$ for $b\in\dom(d_\delta)$. (Once more, by the proof of Proposition \ref{duniprop}, this is well-defined as $d$ and $d_\delta$ coincide on the intersection of their domains.) By induction we show that for every choice of elements $b_1,\ldots,b_n$ from $\dom(d)\cup\dom(d_\delta)$ and for all $x\in\dom(\delta)$ (so that also $xb_1\ldots b_n$ is in $\dom(\delta)$)
\beqn{
\delta(xb_1\ldots b_n)\,-\,\delta(x)b_1\ldots b_n
~~=~~
x\,\Bfam{\,d'(b_1)b_2\ldots b_n~+~\ldots~+~b_1\ldots b_{n-1}d'(b_n)\,}.
}\eeqn
This shows that for every $b$ in the new domain there is a uniqe $b'\in\cB$ satisfying $xb'=\delta(xb)-\delta(x)b$ and that the map $d'\colon b\mapsto b'$ is linear. Clearly, $d'$ is a derivation and $\delta$ is a \nbd{d'}derivation. By Parts \eqref{p2} and \eqref{p3}, $d'$ is closable, if and only if $d_\delta$ or, equivalently, if $\delta$ is closable. In other words, every derivation $d$ that turns a closable ternary derivation $\delta$ of a full Hilbert \nbd{\cB}module into a \nbd{d}derivation admits a unique minimal closed extension $\ol{d'}\supset d_\delta$, and $\delta$ is also a \nbd{\ol{d'}}derivation.
\eob

\bthm\label{dclosethm}
Suppose $\delta$ is the generator of a dynamical system $u$ on a full Hilbert \nbd{\cB}module $E$ and a \nbd{d}derivation for some (by Theorem \ref{genterthm} and Proposition \ref{terdprop}\eqref{p3}, necessarily closable) derivation $d$ of $\cB$. Denote by $d_\vp$ the generator of the dynamical system $\vp$ associated with $u$.

\begin{enumerate}
\item
The unique minimal closed extension $\ol{d'}\supset d_\delta$ of $d$ (see Observation \ref{d'ob}) is the generator of $\vp$, if and only if $d\subset d_\vp$.

\item
If $d\subset d_\vp$, then for $\ol{d}=d_\vp$ it is necessary and sufficient that $\ol{d}\supset d_\delta$.
\end{enumerate}
\ethm

\proof
As in the proof of Theorem \ref{genterthm} we see that $\dom(d_\delta)\subset\dom(d_\vp)$ and that the span of $\AB{\sA(\delta),\sA(\delta)}$ is a dense subspace of entire analytic elements of $\dom(d_\vp)$. Therefore, every subspace $D$ with
\beqn{
\AB{\sA(\delta),\sA(\delta)}
~\subset~
D
~\subset~
\dom(d_\vp)
}\eeqn
is a core for $d_\vp$. In particular, $\dom(d_\delta)$ is a core for $d_\vp$.

\item[(1)]
If $d\nsubset d_\vp$, then $d\subset d'\subset\ol{d'}\nsubset d_\vp$. Conversely, if $d\subset d_\vp$ then also $d'\subset d_\vp$ (because $d_\delta\subset d_\vp$ and, therefore, $\alg^*(\dom(d),\dom(d_\delta))\subset\dom(d_\vp)$), so that $\ol{d'}\subset\ol{d_\vp}=d_\vp$.

\item[(2)]
If $\ol{d}\nsupset d_\delta$, then $\ol{d}\nsupset\ol{d_\delta}=d_\vp$. Conversely, if $\ol{d}\supset d_\delta$ so that $\ol{d}\supset\ol{d_\delta}=d_\vp$, then $d_\vp=\ol{d_\vp}\supset\ol{d}\supset d_\vp$.\qed

\bcor\label{gengencor}
If $\delta$ is the generator of a dynamical system $u$ on a full Hilbert \nbd{\cB}module, then $\delta$ is a \nbd{\ol{d_\delta}}derivation and $\ol{d_\delta}$ is the generator of the \nbd{C^*}dynamical system associated with $u$.
\ecor

In general, a \nbd{d}derivation (even bounded) of a full Hilbert \nbd{\cB}module for some derivation $d$ of $\cB$ need not be a ternary derivation, not even if $d$ is a bounded \nbd{*}derivation.

\bex
The so-called \hl{inner} generalized derivations of a Hilbert \nbd{\cB}module $E$ are the mappings that can be written in the form
\beqn{
\delta(x)
~=~
\alpha x-x\beta
}\eeqn
form some $\alpha\in\sB^a(E)$ and $\beta\in\cB$. From
\bmun{
\delta(x)\AB{y,z}+x\AB{\delta(y),z}+x\AB{y,\delta(z)}
~=~
(\alpha x-x\beta)\AB{y,z}+x\AB{\alpha y-y\beta,z}+x\AB{y,\alpha z-z\beta}
\\
~=~
\alpha x\AB{y,z}-x\beta\AB{y,z}+x\AB{\alpha y,z}-x\AB{y\beta,z}+x\AB{y,\alpha z}-x\AB{y,z\beta}
\\
~=~
\delta(x\AB{y,z})-x\AB{y(\beta+\beta^*),z}+x\AB{y,(\alpha+\alpha^*)z}
}\emun
we see that $\delta$ is a ternary derivation, if and only if $(\beta+\beta^*)\AB{y,z}=\AB{y,(\alpha+\alpha^*)z}$ for all $y,z\in E$. Inserting $yb$ for $y$ and computing $\AB{yb,(\alpha+\alpha^*)z}=b^*\AB{y,(\alpha+\alpha^*)z}$, one may check that $\beta+\beta^*$ must be in the center of $\cB$. Further, the element $\alpha+\alpha^*\in\sB^a(E)$ is given simply as multiplication from the right with the central element $\beta+\beta^*$. Therefore, $\delta$ is a ternary derivation, if and only if the real parts of $\alpha$ and $\beta$ may be removed without changing $\delta$, or, in other words, if $\delta(x)=\alpha x-x\beta$ for skew-adjoint elements $\alpha$ and $\beta$.

Notice, further, that $\delta$ is the generator of the uniformly continuous one-parameter group $u_t(x)=e^{t\alpha}xe^{-t\beta}$ on $E$. It follows that this group is a dynamical system, if and only the groups $e^{t\alpha}$ and $e^{-t\beta}$ are unitary. So, even if $\beta$ is skew-adjoint (so that $d$ is a \nbd{*}derivation and the generator of a \nbd{C^*}dynamical system) $\delta$ does not generate a dynamical system, unless also $\alpha$ is skew-adjoint. On the other hand, if, in this case, $\alpha$ is not skew-adjoint, then $u_t$ is not a \nbd{C_0}group.
\eex

We will see in a moment that the last statement of the preceding example is typical in the sense that, if a \nbd{C_0}group $u$ consists of \nbd{\vp_t}linear maps $u_t$, then $u$ is a dynamical system. But, we think that the following preparatory result inspired very much by Lance \cite[Theorem 3.5]{Lan95} is worth to be stated separately.

\blem\label{philinlem}
Let $E$ be a Hilbert \nbd{\cB}module, let $F$ be a Hilbert \nbd{\cC}module and suppose $u\colon E\rightarrow F$  is a Banach space isometry onto a \nbd{\cC}submodule of $F$. If $u$ is \nbd{\vp}linear for some homomorphism $\vp\colon\cB\rightarrow\cC$ such that $\vp(\cB)\supset\AB{u(F),u(F)}$, then $u$ is a \nbd{\vp}isometry.
\elem

\proof
For $\cC=\cB$, $\vp=\id_\cB$  and surjective $u$ the statement is exactly \cite[Theorem 3.5]{Lan95}. We shall prove the statement exactly along the lines of the proof of \cite[Theorem 3.5]{Lan95} by appealing to \cite[Lemma 3.4]{Lan95} which states
\beqn{
b_1\ge0,b_2\ge0,
~~\text{and}~~
\norm{b_1b}
=
\norm{b_2b}
~\forall~
b\in\cB
~~~~~~\Longrightarrow~~~~~~
b_1
=
b_2.
}\eeqn
First, we compute
\beqn{
\norm{ux}\norm{\vp(b)}
~\ge~
\norm{(ux)\vp(b)}
~=~
\norm{u(xb)}
~=~
\norm{xb}.
}\eeqn
If $0\ne b\in\ol{\cB_E}$ then there exists $x\in E$ such that $xb\ne0$. By Lemma \ref{nlem}, it follows that $\norm{\vp(b)}=\norm{b}$ for all $b\in\ol{\cB_E}$. Next, for all $b\in\cB$ and for all $x\in E$ we have
\bmun{
\norm{\vp(b^*)\AB{ux,ux}\vp(b)}
~=~
\norm{\AB{u(xb),u(xb)}}
~=~
\norm{u(xb)}^2
~=~
\norm{xb}^2
\\
~=~
\norm{b^*\AB{x,x}b}
~=~
\norm{\vp(b^*)\vp(\AB{x,x})\vp(b)},
}\emun
where the last equality follows from $b^*\AB{x,x}b\in\cB_E$ and the first step. In other words, we have $\norm{\sqrt{\AB{ux,ux}}c}=\norm{\sqrt{\vp(\AB{x,x})}c}$ for all elements $c\in\vp(\cB)$. Since by assumption $\AB{ux,ux}\in\vp(\cB)$ so that also $\sqrt{\AB{ux,ux}}\in\vp(\cB)$, it follows by \cite[Lemma 3.4]{Lan95} that $\sqrt{\AB{ux,ux}}=\sqrt{\vp(\AB{x,x})}$, hence, $\AB{ux,ux}=\vp(\AB{x,x})$ and, finally, by polarization $\AB{ux,uy}=\vp(\AB{x,y})$ for all $x,y\in E$. In other words, $u$ is a \nbd{\vp}isometry.\qed

\bcor\label{philincor}
Every \nbd{\vp}linear, isometric Banach space isomorphism between full Hilbert modules with surjective $\vp$ is necessarily a \nbd{\vp}unitary.
\ecor

\brem
We do not know, whether the (necessary) condition $\vp(\cB)\supset\AB{u(F),u(F)}$ in Lemma \ref{philinlem} (and the corresponding condition $\vp$ be surjective of Corollary \ref{philincor}) does not, possibly, follow from the remaining hypothesis.
\erem

\bthm\label{gengenthm}
Suppose that $d$ is a \nbd{*}derivation that is the generator of a \nbd{C_0}group $\vp$ on the \nbd{C^*}al\-ge\-bra $\cB$, and suppose that $\delta$ is a \nbd{d}derivation that is the generator of a \nbd{C_0}group $u$ on the full Hilbert \nbd{\cB}module $E$. Then $u$ is a dynamical system on $E$ and $\vp$ is the \nbd{C^*}dynamical system associated with $u$. Of course, $\cB_E$ is a core for $d$ and $\delta$ is a ternary derivation and a \nbd{d}derivation.
\ethm

\proof
For all $x\in\sA(\delta)$ and all $b\in\sA(d)$ as in the proof of Theorem \ref{genterthm} one shows that also $xb\in\sA(\delta)$ and that
\beqn{
u_t(xb)
~=~
u_t(x)\vp_t(b).
}\eeqn
In exactly the same way one shows that $\vp_t$ (of course, a \nbd{*}map) is multiplicative. In other words, $\vp_t$ is an automorphism of $\cB$ and $u_t$ is a surjective and right \nbd{\vp_t}linear Banach space isometry. By Corollary \ref{philincor}, $u_t$ is a \nbd{\vp_t}unitary. In other words, $u$ is a dynamical system and $\vp$ is the \nbd{C^*}dynamical system associated with it.\qed

\lf
For the sake of clarity we summarize the criteria provided by Theorem \ref{genterthm}, Corollary \ref{gengencor}, and Theorem \ref{gengenthm}. Without the obvious proof, we add a fourth criterion based on the observation (as explained in Remark \ref{linkrem}) that a dynamical system on $E$ extends to a \nbd{C^*}dynamical system on the linking algebra.

\bthm\label{sumthm}
Let $\delta$ be the generator of a \nbd{C_0}group $u$ on a full Hilbert \nbd{\cB}module. Then the following statements are equivalent:
\begin{enumerate}
\item
$u$ is a dynamical system.

\item
$\delta$ is a ternary derivation.

\item
There exists a \nbd{*}derivation $d$ that is the generator of a \nbd{C_0}group on $\cB$ (necessarily a \nbd{C^*}dy\-nam\-ical system) such that $\delta$ is a \nbd{d}derivation.

\item
$\delta$ extends to the generator of a \nbd{C^*}dynamical system on the linking algebra of the form $\Delta=\tMatrix{d&\delta^*\\\delta&D}$ with $\delta(x^*)^*:=\delta(x)^*$ and $d$ and $D$ being generators of \nbd{C^*}dynamical systems on $\cB$ and $\sK(E)$, respectively.
\end{enumerate}
\ethm

\brem
In all criteria where we make explicit reference to a derivation $d$ of the \it{corner} $\cB$, we assume that both $\delta$ and $d$ are generators of \nbd{C_0}groups. We leave open the very interesting question whether the algebraic conditions alone might already be sufficient to conclude from one, $\delta$ or $d$, being a generator, that also the other is a generator.
\erem

\setlength{\baselineskip}{2.5ex}


\newcommand{\Swap}[2]{#2#1}\newcommand{\Sort}[1]{}
\providecommand{\bysame}{\leavevmode\hbox to3em{\hrulefill}\thinspace}
\providecommand{\MR}{\relax\ifhmode\unskip\space\fi MR }
\providecommand{\MRhref}[2]{%
  \href{http://www.ams.org/mathscinet-getitem?mr=#1}{#2}
}
\providecommand{\href}[2]{#2}

\noindent
Gholamreza Abbaspour Tabadkan: {\small\itshape Department of Mathematics}, {\small\itshape Ferdowsi University}, {\small\itshape P.O.\ Box 1159, Mashhad, Iran},
{\small\itshape E-mail: \tt{tabadkan@math.um.ac.ir}}

\lf\noindent
Michael Skeide: {\small\itshape Dipartimento S.E.G.e S.}, {\small\itshape Università degli Studi del Molise}, {\small\itshape Via de Sanctis}, {\small\itshape 86100 Campobasso, Italy}, {\small{\itshape E-mail: \tt{skeide@math.tu-cottbus.de}}},
\\
{\small{\itshape Homepage: \tt{http://www.math.tu-cottbus.de/INSTITUT/lswas/\_skeide.html}}}


\end{document}